\documentclass[12pt]{amsart}

\usepackage{amssymb}
\usepackage[bbgreekl]{mathbbol}
\usepackage{graphicx}
\usepackage{hyperref}
\usepackage{ifthen}
\usepackage{mathrsfs}
\usepackage{pict2e}
\usepackage{xargs}
\usepackage{xspace}

\DeclareSymbolFontAlphabet{\amsmathbb}{AMSb}

\newcommand{\definedterm}[1]{\emph{#1}}

\newcommand{\action}{\curvearrowright}
\newcommand{\Bairespace}[1][]{
  \ifthenelse{\equal{#1}{}}{\functions{\N}{\N}}{\functions{#1}{\N}}
}
\newcommand{\Bairetree}[1][]{
  \ifthenelse{\equal{#1}{}}{\functions{<\N}{\N}}{\functions{#1}{\N}}
}

\newcommand{\bbG}{\amsmathbb{G}}
\newcommand{\bbs}[1][]{\mathbb{s}_{#1}}

\newcommand{\calN}{\mathcal{N}}

\newcommand{\Cantorspace}[1][]{
  \ifthenelse{\equal{#1}{}}{\functions{\N}{2}}{\functions{#1}{2}}
}
\newcommand{\Cantortree}[1][]{
  \ifthenelse{\equal{#1}{}}{\functions{<\N}{2}}{\functions{#1}{2}}
}

\newcommand{\composition}{\circ}
\newcommandx{\concatenation}[2][1 = undefined, 2 = undefined]{
  \ifthenelse{\equal{#1}{undefined}}{{}\smallfrown}{
    \ifthenelse
      {\equal{#2}{undefined}}
      {\smallfrown_{#1}}
      {\bigoplus_{#1} #2}
  }
}
\newcommand{\continuum}{2^{\aleph_0}}
\newcommand{\emptysequence}{\emptyset}
\newcommand{\equivalenceclass}[2]{[#1]_{#2}}
\newcommand{\extendedby}{\sqsubseteq}
\newcommand{\extensions}[1]{\calN_{#1}}
\newcommand{\from}{\colon}
\newcommandx{\functions}[3][3 =]{
  \ifthenelse{\equal{#3}{}}{#2^{#1}}{#2^{#1}_{#3}}
}
\newcommand{\Gdelta}{$G_\delta$\xspace}
\newcommand{\goesto}{\rightarrow}
\newcommand{\Gzero}{\bbG_0}
\newcommand{\image}[2]{#1(#2)}
\newcommandx{\intersection}[2][1 =, 2 =]{
  \ifthenelse{\equal{#1}{}}{\cap}{
    \ifthenelse{\equal{#2}{}}{\bigcap #1}{{\bigcap_{#1} #2}}
  }
}
\newcommand{\into}{\hookrightarrow}

\newcommand{\mathand}{\text{ and }}
\newcommand{\metric}[2][]{
  \ifthenelse{\equal{#1}{}}{d_{#2}}{d_{#1}^{#2}}
}
\newcommand{\N}{\amsmathbb{N}}

\newcommand{\onto}{\twoheadrightarrow}

\newcommand{\pair}[2]{(#1, #2)}
\newcommand{\preimage}[2]{#1^{-1}(#2)}
\newcommandx{\product}[2][1 =, 2 =]{
  \ifthenelse{\equal{#1}{}}{\times}{
    \ifthenelse{\equal{#2}{}}{\prod #1}{{\prod_{#1} #2}}
  }
}
\newcommand{\projection}[1][]{
  \ifthenelse{\equal{#1}{}}{\mathrm{proj}}{\mathrm{proj}_{#1}}
}
\renewcommand{\restriction}[2]{#1 \upharpoonright #2}
\newcommandx{\sequence}[2][2 = undefined]{
  \ifthenelse{\equal{#2}{undefined}}{(#1)}{
    (#1)_{#2}
  }
}
\newcommandx{\set}[2][2 = undefined]{
  \ifthenelse{\equal{#2}{undefined}}{\{ #1 \}}{
    \{ #1 \suchthat #2 \}
  }
}
\newcommand{\setcomplement}[1]{\twiddle #1}
\newcommand{\strictlyextendedby}{\sqsubset}
\newcommand{\suchthat}{\mid}

\newcommand{\triple}[3]{(#1, #2, #3)}
\newcommand{\twiddle}
  {\raisebox{1.5pt}{\scalebox{.75}{$\mathord{\sim}$}}}
\newcommandx{\union}[2][1 =, 2 =]{
  \ifthenelse{\equal{#1}{}}{\cup}{
    \ifthenelse{\equal{#2}{}}{\bigcup #1}{{\bigcup_{#1} #2}}
  }
}

\newcommand{\Akin}{Akin\xspace}
\newcommand{\Baire}{Baire\xspace}
\newcommand{\Blanchard}{Blan\-chard\xspace}
\newcommand{\Borel}{Bor\-el\xspace}
\newcommand{\Cantor}{Can\-tor\xspace}

\newcommand{\Hausdorff}{Haus\-dorff\xspace}
\newcommand{\Huang}{Hu\-ang\xspace}
\newcommand{\Kechris}{Kech\-ris\xspace}
\newcommand{\Kuratowski}{Kur\-at\-ow\-ski\xspace}
\newcommand{\Li}{Li\xspace}
\newcommand{\Mycielski}{My\-ciel\-ski\xspace}
\newcommand{\Polish}{Po\-lish\xspace}

\newcommand{\Silver}{Sil\-ver\xspace}
\newcommand{\Snoha}{Sno\-ha\xspace}
\newcommand{\Solecki}{Sol\-eck\-i\xspace}
\newcommand{\Souslin}{Sous\-lin\xspace}
\newcommand{\Todorcevic}{To\-dor\-cev\-ic\xspace}
\newcommand{\Ulam}{U\-lam\xspace}
\newcommand{\Yorke}{Yorke\xspace}

\newenvironment{lemmaproof}{
  
  \begin{proof}
}{\end{proof}}

\newenvironment{propositionproof}{
  
  \begin{proof}
}{\end{proof}}

\newenvironment{theoremproof}{
  
  \begin{proof}
}{\end{proof}}

\newtheorem{introtheorem}{Theorem}

\newtheorem{lemma}{Lemma}[section]
\newtheorem{proposition}[lemma]{Proposition}
\newtheorem{theorem}[lemma]{Theorem}

\begin{document}


\begin{abstract}
  We show that \Li--\Yorke chaos ensures the existence of a scrambled
  Cantor set.
\end{abstract}

\author[S. Geschke]{Stefan Geschke}
\address{
  Stefan Geschke \\
  Universit\"{a}t Hamburg \\
  Department of Mathematics \\
  Bundesstrasse 55 \\
  20146 Hamburg \\
  Germany
 }
\email{stefan.geschke@uni-hamburg.de}
\urladdr{
  \url{https://www.math.uni-hamburg.de/home/geschke/}
}
\thanks{The first author was partially supported by DAAD project
  grant 57156702.}
  
\author{Jan Greb\'ik}
\address{
Jan Greb\'ik\\
Mathematics Institute\\
University of Warwick\\
Coventry CV4 7AL, UK
}
\email{jan.grebik@warwick.ac.uk}
\urladdr{
  \url{http://homepages.warwick.ac.uk/staff/Jan.Grebik/}
 }
 \thanks{The second author was partially supported by DAAD project
  grant 15-13.}

\author[B.D. Miller]{Benjamin D. Miller}
\address{
  Benjamin D. Miller \\
  Universit\"{a}t Wien \\
  Department of Mathematics \\
  Oskar Morgenstern Platz 1 \\
  1090 Wien \\
  Austria
 }
\email{benjamin.miller@univie.ac.at}
\urladdr{
  \url{https://homepage.univie.ac.at/benjamin.miller/}
}
\thanks{The third author was partially supported by FWF stand-alone
  project grants P28153 and P29999.}



\title{Scrambled Cantor sets}

\maketitle

\section*{Introduction}

A \definedterm{dynamical system} is a pair $\pair{X}{f}$, where $X$ is
a metric space and $f \from X \to X$ is a continuous function. Given
such a system, we say that points $x, y \in X$ are \definedterm
{proximal} if $\liminf_{n \to \infty} \metric{X}(f^n(x), f^n(y)) = 0$, and
\definedterm{asymptotic} if $\limsup_{n \to \infty} \metric{X}(f^n(x),
f^n(y)) = 0$. The pair $\pair{x}{y}$ is \definedterm{\Li--\Yorke} if $x$ and
$y$ are proximal but not asymptotic, a set $Y \subseteq X$ is
\definedterm{scrambled} if $\pair{x}{y}$ is \Li--\Yorke for all distinct $x, y
\in Y$, and the system $\pair{X}{f}$ is \definedterm{\Li--\Yorke chaotic} if
there is an uncountable scrambled set $Y \subseteq X$. In \cite{LY}, \Li
and \Yorke showed that every dynamical system on the unit interval with
a point of period three is \Li--\Yorke chaotic.

The scrambled set constructed in \cite{LY} is indexed by an interval on
the real line, and therefore has cardinality $\continuum$. Moreover,
subsequent constructions of uncountable scrambled sets in the
literature typically gave rise to sets of cardinality $\continuum$, or even
\definedterm{Cantor sets}, i.e., homeomorphic copies of the Cantor
space $\Cantorspace$. One example is the construction, in \cite
{BGKM}, of uncountable scrambled sets in dynamical systems of
positive topological entropy.

A metric space is \definedterm{\Polish} if it is complete and separable,
and a dynamical system $\pair{X}{f}$ is \definedterm{\Polish} if $X$ is
\Polish. At the end of \cite[\S3]{BHS}, \Blanchard, \Huang, and \Snoha
asked whether every \Li--\Yorke chaotic \Polish dynamical system
admits a scrambled \Cantor set. Here we utilize the descriptive set
theory of definable graphs to obtain the following answer to their
question:

\begin{introtheorem} \label{intro:uncountable}
  Suppose that $\pair{X}{f}$ is a \Li--\Yorke chaotic \Polish dynamical
  system. Then there is a scrambled Cantor set $C \subseteq X$.
\end{introtheorem}

In \S\ref{gzero}, we establish an analog of the
\Kechris--\Solecki--\Todorcevic characterization of the existence of
$\aleph_0$-colorings (see \cite[Theorem 6.3]{KST}) within cliques. In
\S\ref{perfect}, we use this to establish a similar analog of \Silver's
perfect set theorem for co-analytic equivalence relations (see
\cite{Silver}). In \S\ref{chaos}, we use the latter to establish Theorem
\ref{intro:uncountable}. And in \S\ref{generalizations}, we discuss
several generalizations.

\section{Colorings in cliques} \label{gzero}

Endow $\N$ with the discrete topology, and $\Bairespace$ with the
corresponding product topology. A topological space is \definedterm
{analytic} if it is a continuous image of a closed subset of
$\Bairespace$, and \definedterm{\Polish} if it is second countable
and admits a compatible complete metric. A subset of a topological
space is \definedterm{\Borel} if it is in the smallest $\sigma$-algebra
containing the open sets, and \definedterm{co-analytic} if its
complement is analytic. Every non-empty \Polish space is a
continuous image of $\Bairespace$ (see, for example, \cite[Theorem
7.9]{Kechris}), thus so too is every non-empty analytic space, and a
subset of an analytic \Hausdorff space is \Borel if and only if it is
analytic and co-analytic (see, for example, the proof of \cite[Theorem
14.11]{Kechris}).

A \definedterm{digraph} on a set $X$ is an irreflexive binary relation
$G$ on $X$. A set $Y \subseteq X$ is \definedterm{$G$-independent}
if $\restriction{G}{Y} = \emptyset$. An \definedterm{$I$-coloring} of $G$
is a function $c \from X \to I$ such that $c(x) \neq c(y)$ for all $\pair{x}
{y} \in G$, or equivalently, such that $\preimage{c}{\set{i}}$ is
$G$-independent for all $i \in I$. A \definedterm{homomorphism} from a
binary relation $R$ on a set $X$ to a binary relation $S$ on a set $Y$ is
a function $\phi \from X \to Y$ for which $\image{(\phi \times \phi)}{R}
\subseteq S$. We say that a set $Y \subseteq X$ is an \definedterm
{$R$-clique} if $x \mathrel{R} y$ for all distinct $x, y \in Y$.

We use $\functions{<\N}{X}$ to denote $\union[n \in \N][\functions{n}
{X}]$, $\sequence{i}$ to denote the singleton sequence with value $i$,
and $\extendedby$ to denote extension. Following standard practice,
we also use $\extensions{s}$ to denote $\set{b \in \Bairespace}[s
\extendedby b]$ or $\set{c \in \Cantorspace}[s \extendedby c]$ (with
the context determining which of the two we have in mind). Fix
sequences $\bbs[n] \in \Cantorspace[n]$ such that $\forall s \in
\Cantorspace \exists n \in \N \ s \extendedby \bbs[n]$, and let $\Gzero$
denote the digraph on $\Cantorspace$ given by $\Gzero = \set
{\sequence{\bbs[n] \concatenation \sequence{i} \concatenation c}[i < 2]}
[c \in \Cantorspace \mathand n \in \N]$. A subset of a topological space
is \Gdelta if it is a countable intersection of open sets.

\begin{theorem}
  \label{gzero:dichotomy}
  Suppose that $X$ is a \Hausdorff space, $G$ is an analytic digraph on
  $X$, and $R$ is a reflexive \Gdelta binary relation on $X$. Then at
  least one of the following holds:
  \begin{enumerate}
    \item For every $R$-clique $Y \subseteq X$, there is an $\N$-coloring
      of $\restriction{G}{Y}$.
    \item There is a continuous homomorphism $\phi \from \Cantorspace
      \to X$ from $\Gzero$ to $G$ for which $\image{\phi}{\Cantorspace}$
        is an $R$-clique.
  \end{enumerate}
\end{theorem}

\begin{theoremproof}
  Suppose that there is an $R$-clique $Y \subseteq X$ for which there
  is no $\N$-coloring of $\restriction{G}{Y}$. Then $G \neq \emptyset$,
  so there are continuous surjections $\phi_G \from \Bairespace \onto
  G$ and $\phi_X \from \Bairespace \onto \union[i < 2][\image{\projection
  [i]}{G}]$. Fix a decreasing sequence $\sequence{R_n}[n \in \N]$ of
  open subsets of $X \times X$ whose intersection is $R$.
   
  We will define a decreasing sequence $\sequence{Y^\alpha}[\alpha <
  \omega_1]$ of subsets of $Y$, off of which there are $\N$-colorings
  of $\restriction{G}{Y}$. We begin by setting $Y^0 = Y$. For all limit
  ordinals $\lambda < \omega_1$, we set $Y^\lambda = \intersection
  [\alpha < \lambda][Y^\alpha]$. To describe the construction at
  successor ordinals, we require several preliminaries.
  
  An \definedterm{approximation} is a triple of the form $a = \triple{n^a}
  {\phi^a}{\sequence{\psi_n^a}[n < n^a]}$, where $n^a \in \N$, $\phi^a
  \from \Cantorspace[n^a] \to \Bairetree$, $\psi_n^a \from \Cantorspace
  [n^a - (n+1)] \to \functions{n^a}{\N}$ for all $n < n^a$, and $\image
  {\phi_X}{\extensions{\phi^a(s)}} \times \image{\phi_X}{\extensions
  {\phi^a(t)}} \subseteq R_{n^a}$ for all distinct $s, t \in \Cantorspace
  [n^a]$. A \definedterm{one-step extension} of an approximation $a$ is
  an approximation $b$ such that:
  \begin{enumerate}
    \renewcommand{\theenumi}{\alph{enumi}}
    \setlength{\itemindent}{-6pt}
    \item $n^b = n^a + 1$.
    \item $\forall s \in \Cantorspace[n^a] \forall t \in \Cantorspace[n^b]
      \ (s \strictlyextendedby t \implies \phi^a(s) \strictlyextendedby
        \phi^b(t))$.
    \item $\forall n < n^a \forall s \in \Cantorspace[n^a - (n + 1)] \forall t
      \in \Cantorspace[n^b - (n + 1)] \ (s \strictlyextendedby t \implies
        \psi^a_n(s) \strictlyextendedby \psi^b_n(t))$.
  \end{enumerate}
  Similarly, a \definedterm{configuration} is a triple of the form $\gamma
  = \triple{n^\gamma}{\phi^\gamma}{\sequence{\psi_n^\gamma}[n <
  n^\gamma]}$, where $n^\gamma \in \N$, $\phi^\gamma \from
  \Cantorspace[n^\gamma] \to \Bairespace$, $\psi_n^\gamma \from
  \Cantorspace[n^\gamma - (n+1)] \to \Bairespace$ for all $n <
  n^\gamma$, and $(\phi_G \composition \psi_n^\gamma)(t) =
  \sequence{(\phi_X \composition \phi^\gamma)(\bbs[n] \concatenation
  \sequence{i} \concatenation t)}[i < 2]$ for all $n < n^\gamma$ and $t
  \in \Cantorspace[n^\gamma - (n + 1)]$. We say that $\gamma$ is
  \definedterm{compatible} with a set $Y' \subseteq Y$ if $\image
  {(\phi_X \composition \phi^\gamma)}{\Cantorspace[n^\gamma]}
  \subseteq Y'$, and \definedterm{compatible} with $a$ if:
  \begin{enumerate}
    \renewcommand{\theenumi}{\roman{enumi}}
    \item $n^a = n^\gamma$.
    \item $\forall t \in \Cantorspace[n^a] \ \phi^a(t) \extendedby
      \phi^\gamma(t)$.
    \item $\forall n < n^a \forall t \in \Cantorspace[n^a - (n + 1)] \ \psi^a_n
      (t) \extendedby \psi^\gamma_n(t)$.
  \end{enumerate}
  An approximation $a$ is \definedterm{$Y'$-terminal} if no configuration
  is compatible with both $Y'$ and a one-step extension of $a$. Let $A
  (a, Y')$ denote the set of points of the form $(\phi_X \composition
  \phi^\gamma)(\bbs[n^a])$, where $\gamma$ varies over all
  configurations compatible with $a$ and $Y'$.
  
  \begin{lemma} \label{chromaticnumber:terminalimpliesindependent}
    Suppose that $Y' \subseteq Y$ and $a$ is a $Y'$-terminal
    approximation. Then $A(a, Y')$ is $G$-independent.
  \end{lemma}
  
  \begin{lemmaproof}
    Suppose, towards a contradiction, that there are configurations
    $\gamma_0$ and $\gamma_1$, both compatible with $a$ and $Y'$,
    with the property that $\sequence{(\phi_X \composition
    \phi^{\gamma_i})(\bbs[n^a])}[i < 2] \in G$. Fix a sequence $d \in
    \Bairespace$ with the property that $\phi_G(d) = \sequence{(\phi_X
    \composition \phi^{\gamma_i})(\bbs[n^a])}[i < 2]$, and let $\gamma$
    be the configuration given by $n^\gamma = n^a + 1$, $\phi^\gamma
    (t \concatenation \sequence{i}) = \phi^{\gamma_i}(t)$ for all $i < 2$
    and $t \in \Cantorspace[n^a]$, $\psi_n^\gamma(t \concatenation
    \sequence{i}) = \psi_n^{\gamma_i}(t)$ for all $i < 2$, $n < n^a$, and
    $t \in \Cantorspace[n^a - (n + 1)]$, and $\psi_{n^a}^\gamma
    (\emptysequence) = d$. Then $\gamma$ is compatible with a
    one-step extension of $a$, contradicting the fact that $a$ is
    $Y'$-terminal.
  \end{lemmaproof}
  
  Let $Y^{\alpha + 1}$ be the difference of $Y^\alpha$ and the union of
  the sets of the form $A(a, Y^\alpha)$, where $a$ varies over all
  $Y^\alpha$-terminal approximations.
  
  \begin{lemma}
    \label{chromaticnumber:nonterminalapproximationextension}
    Suppose that $\alpha < \omega_1$ and $a$ is a non-$Y^{\alpha +
    1}$-terminal approximation. Then $a$ has a non-$Y^\alpha$-terminal
    one-step extension.
  \end{lemma}
  
  \begin{lemmaproof}
    Fix a one-step extension $b$ of $a$ for which there is a configuration
    $\gamma$ compatible with $b$ and $Y^{\alpha + 1}$. Then $(\phi_X
    \composition \phi^\gamma)(\bbs[n^b]) \in Y^{\alpha + 1}$, so $b$ is
    not $Y^\alpha$-terminal. 
  \end{lemmaproof}
  
  Fix $\alpha < \omega_1$ such that the families of $Y^\alpha$- and
  $Y^{\alpha + 1}$-terminal approximations coincide, and let $a_0$ be
  the unique approximation for which $n^{a_0} = 0$ and $\phi^{a_0}
  (\emptysequence) = \emptysequence$. As $A(a_0, Y') = Y'
  \intersection \union[i < 2][\image{\projection[i]}{G}]$ for all $Y'
  \subseteq Y$,  we can assume that $a_0$ is not $Y^\alpha$-terminal,
  since otherwise there is an $\N$-coloring of $\restriction{G}{Y}$.
  
  By recursively applying Lemma \ref
  {chromaticnumber:nonterminalapproximationextension}, we obtain
  non-$Y^\alpha$-terminal one-step extensions $a_{n+1}$ of $a_n$ for
  all $n \in \N$. Define $\phi', \psi_n \from \Cantorspace \to \Bairespace$
  by $\phi'(c) = \union[n \in \N][\phi^{a_n}(\restriction{c}{n})]$ and $\psi_n
  (c) = \union[m > n][\psi_n^{a_m}(\restriction{c}{(m - (n + 1))})]$ for all
  $c \in \Cantorspace$ and $n \in \N$. Clearly these functions are
  continuous.
  
  To see that the function $\phi = \phi_X \composition \phi'$ is a
  homomorphism from $\Gzero$ to $G$, we will show the stronger fact
  that if $c \in \Cantorspace$ and $n \in \N$, then $(\phi_G \composition
  \psi_n)(c) = \sequence{(\phi_X \composition \phi')(\bbs[n]
  \concatenation \sequence{i} \concatenation c)}[i < 2]$. It is sufficient to
  show that if $U$ is an open neighborhood of $\sequence{(\phi_X
  \composition \phi')(\bbs[n] \concatenation \sequence{i} \concatenation
  c)}[i < 2]$ and $V$ is an open neighborhood of $(\phi_G \composition
  \psi_n)(c)$, then $U \intersection V \neq \emptyset$. Towards this end,
  fix $m > n$ such that $\product[i < 2][\image{\phi_X}{\extensions
  {\phi^{a_m}(\bbs[n] \concatenation \sequence{i} \concatenation s)}}]
  \subseteq U$ and $\image{\phi_G}{\extensions{\psi_n^{a_m}(s)}}
  \subseteq V$, where $s = \restriction{c}{(m - (n + 1))}$. As $a_m$ is
  not $Y^\alpha$-terminal, there is a configuration $\gamma$ compatible
  with $a_m$, so $\sequence{(\phi_X \composition \phi^\gamma)(\bbs[n]
  \concatenation \sequence{i} \concatenation s)}[i < 2] \in U$ and
  $(\phi_G \composition \psi_n^\gamma)(s) \in V$, thus $U \intersection
  V \neq \emptyset$.
  
  To see that $\image{\phi}{\Cantorspace}$ is an $R$-clique, observe
  that if $c, d \in \Cantorspace$ are distinct and $n \in \N$ is sufficiently
  large that $\restriction{c}{n} \neq \restriction{d}{n}$, then $\phi(c) \in
  \image{\phi_X}{\extensions{\phi^{a_n}(\restriction{c}{n})}}$ and $\phi
  (d) \in \image{\phi_X}{\extensions{\phi^{a_n}(\restriction{d}{n})}}$, so
  $\phi(c) \mathrel{R_n} \phi(d)$.
\end{theoremproof}

\section{Separability in cliques} \label{perfect}

The following well-known fact rules out the existence of a
\Baire-measurable $\N$-coloring of $\Gzero$:

\begin{proposition} \label{perfect:meager}
  Suppose that $B \subseteq \Cantorspace$ is a non-meager set with
  the \Baire property. Then $B$ is not $\Gzero$-independent.
\end{proposition}

\begin{propositionproof}
  Fix a sequence $s \in \Cantortree$ for which $B$ is comeager in
  $\extensions{s}$ (see, for example, \cite[Proposition 8.26]{Kechris}).
  Then there exists $n \in \N$ for which $s \extendedby \bbs[n]$. Let
  $\iota$ be the involution of $\extensions{\bbs[n]}$ sending $\bbs[n]
  \concatenation \sequence{0} \concatenation c$ to $\bbs[n]
  \concatenation \sequence{1} \concatenation c$ for all $c \in
  \Cantorspace$. As $\iota$ is a homeomorphism, it follows that $B
  \intersection \image{\iota}{B}$ is comeager in $\extensions{\bbs[n]}$
  (see, for example, \cite[Exercise 8.45]{Kechris}), so $B \intersection
  \image{\iota}{B} \intersection \extensions{\bbs[n] \concatenation
  \sequence{0}} \neq \emptyset$. As $\pair{c}{\iota(c)} \in \restriction
  {\Gzero}{B}$ for all $c \in B \intersection \image{\iota}{B} \intersection
  \extensions{\bbs[n] \concatenation \sequence{0}}$, it follows that $B$
  is not $\Gzero$-independent.
\end{propositionproof}

The following corollary is also well known:

\begin{proposition} \label{perfect:meager:product}
  Suppose that $E$ is a non-meager equivalence relation on
  $\Cantorspace$ with the \Baire property. Then $E$ is not disjoint from
  $\Gzero$.
\end{proposition}

\begin{propositionproof}
  By the \Kuratowski-\Ulam theorem (see, for example, \cite[Theorem
  8.41]{Kechris}), there is a sequence $c \in \Cantorspace$ for which
  $\equivalenceclass{c}{E}$ has the \Baire property and is not meager,
  so Proposition \ref{perfect:meager} ensures that $\equivalenceclass{c}
  {E}$ is not $\Gzero$-independent.
\end{propositionproof}

A \definedterm{partial transversal} of an equivalence relation $E$ on a
set $X$ is a set $Y \subseteq X$ that does not contain distinct
$E$-related points.

\begin{theorem} \label{perfect:dichotomy}
  Suppose that $X$ is a \Hausdorff space, $E$ is a co-analytic
  equivalence relation on $X$, and $R$ is a reflexive \Gdelta binary
  relation on $X$ for which there is an $R$-clique $Y \subseteq X$
  intersecting uncountably-many $E$-classes. Then there is a \Cantor
  set $C \subseteq X$ that is both a partial transversal of $E$ and an
  $R$-clique.
\end{theorem}

\begin{theoremproof}
  Set $G = \setcomplement{E}$, appeal to Theorem \ref
  {gzero:dichotomy} to obtain a continuous homomorphism $\phi \from
  \Cantorspace \to X$ from $\Gzero$ to $G$ for which $\image{\phi}
  {\Cantorspace}$ is an $R$-clique, and let $F$ be the pullback of $E$
  through $\phi$. As $\Gzero \intersection F = \emptyset$, Proposition
  \ref{perfect:meager:product} implies that $F$ is meager, thus \Mycielski's
  Theorem (see, for example, \cite[Theorem 19.1]{Kechris}) yields a
  continuous injection $\psi \from \Cantorspace \into \Cantorspace$ sending
  distinct elements of $\Cantorspace$ to $F$-inequivalent elements of
  $\Cantorspace$, in which case the set $C = \image{(\phi \composition
  \psi)}{\Cantorspace}$ is as desired.
\end{theoremproof}

\section{Li-Yorke chaos} \label{chaos}

We say that a dynamical system $\pair{X}{f}$ is \definedterm{analytic}
if $X$ is analytic. As every \Polish dynamical system is analytic, Theorem
\ref{intro:uncountable} is a consequence of the following result:

\begin{theorem} \label{chaos:uncountable}
  Suppose that $\pair{X}{f}$ is a \Li--\Yorke chaotic analytic dynamical
  system. Then there is a scrambled Cantor set $C \subseteq X$.
\end{theorem}

\begin{theoremproof}
  Note first that the set
  \begin{align*}
    R
      & = \set{\pair{x}{y} \in X \times X}[x \mathand y \text{ are proximal}]
        \\
      & = \intersection[\epsilon > 0][{\intersection[n \in \N][{\union[m \ge n]
         [{\set{\pair{x}{y} \in X \times X}[\metric{X}(f^m(x), f^m(y)) <
           \epsilon]}]}]}]
  \end{align*}
  is \Gdelta, and the equivalence relation
  \begin{align*}
    E
      & = \set{\pair{x}{y} \in X \times X}[x \mathand y \text{ are
        asymptotic}] \\
      & = \intersection[\epsilon > 0][{\union[n \in \N][{\intersection[m \ge n]
         [{\set{\pair{x}{y} \in X \times X}[\metric{X}(f^m(x), f^m(y)) \le
           \epsilon]}]}]}]
  \end{align*}
  is \Borel. As the fact that $\pair{X}{f}$ is \Li--\Yorke chaotic yields an
  $R$-clique intersecting uncountably-many $E$-classes, Theorem \ref
  {perfect:dichotomy} yields a scrambled \Cantor set.
\end{theoremproof}

\section{Further generalizations} \label{generalizations}

Suppose that $S$ is a semigroup, $d_S$ is a metric on $S$, $X$ is a
metric space, and $S \action X$ is an action by continuous functions.
The notions of proximal and asymptotic, and therefore of \Li--\Yorke
pair, scrambled, and \Li--\Yorke chaotic, generalize naturally to such
actions, as does Theorem \ref{chaos:uncountable} and its proof.

Given $\delta \ge 0$ and a dynamical system $\pair{X}{f}$, we say that
points $x, y \in X$ are \definedterm{$\delta$-proximal} if $\liminf_{n
\goesto \infty} \metric{X}(f^n(x), f^n(y)) \le \delta$. The above argument
also yields the generalization of Theorem \ref{chaos:uncountable} to the
analog of \Li--\Yorke chaos where proximality is replaced with
$\delta$-proximality.

Given $\epsilon \ge 0$ and a dynamical system $\pair{X}{f}$, we say
that points $x, y \in X$ are \definedterm{$\epsilon$-asymptotic} if
$\limsup_{n \goesto \infty} \metric{X}(f^n(x),f^n(y)) \le \epsilon$, which
gives rise to an analogous notion of \definedterm{$\epsilon$-scrambled}.
By replacing Theorem \ref{perfect:dichotomy} with its generalization
from equivalence relations to extended-valued quasi-metrics in the
proof of Theorem \ref{chaos:uncountable}, one can establish the
generalization of the latter in which the uncountable set is
$\epsilon$-scrambled and the \Cantor set is  $(\epsilon / 2)$-scrambled.

However, when $\epsilon > 0$, the strengthening in which $\epsilon / 2$
is replaced with any $\epsilon' < \epsilon$ follows from the analog of
Theorem \ref{chaos:uncountable} in which asymptoticity is replaced with
the strengthening of $\epsilon$-asymptoticity where one requires that
$\limsup_{n \goesto \infty} \metric{X}(f^n(x),f^n(y)) < \epsilon$, which
was originally established by \Blanchard, \Huang, and \Snoha (see \cite
[Theorem 16]{BHS}). As the corresponding analog of \Li--\Yorke pair is
a \Gdelta condition, their result follows from the special case of
Theorem \ref{perfect:dichotomy} where $E$ is equality, which is far
simpler to establish (see \cite[Remark 1.14]{Shelah}).

At the end of \cite[\S6]{Akin}, \Akin noted that the analog of Theorem
\ref{chaos:uncountable} for complete perfect metric spaces is open.
While our approach does not fully resolve this problem beyond
the separable case, it does generalize to show that if $\kappa$ is an
infinite cardinal, $X$ has a dense set of cardinality $\kappa$, and there
is a scrambled set of cardinality strictly greater than $\kappa$, then
there is a scrambled \Cantor set. To see this, endow $\kappa$ with the
discrete topology, and $\functions{\N}{\kappa}$ with the corresponding
product topology. We say that a topological space is \definedterm
{$\kappa$-\Souslin} if it is a continuous image of a closed subset of
$\functions{\N}{\kappa}$. The proof that every \Polish space is analytic
easily adapts to show that every complete metric space with a dense
subset of cardinality $\kappa$ is $\kappa$-\Souslin, and the proof of
Theorem \ref{gzero:dichotomy} easily adapts to show the analogous
result in which the digraph $G$ is merely $\kappa$-\Souslin and the
$\N$-coloring in condition (1) is replaced with a $\kappa$-coloring. The
proof of Theorem \ref{perfect:dichotomy} therefore adapts to show the
analogous result in which $E$ is co-$\kappa$-\Souslin, the pullback of
$E$ through every continuous function $\phi \from \Cantorspace \to X$
has the \Baire property, and $Y$ intersects more than $\kappa$
classes. But this can be plugged into the proof of Theorem \ref
{chaos:uncountable} to obtain the desired result.

\bibliographystyle{amsalpha}
\bibliography{bibliography}

\end{document}